\documentclass[12pt]{elsarticle}

\usepackage{amssymb}
\journal{}
\usepackage{amssymb}

\usepackage{rotating}

\usepackage{natbib}

\usepackage{multirow}

\usepackage{color}

\usepackage{soul} %VICTOR

\newcommand{\bfg}[1]{\mbox{\boldmath $#1$\unboldmath}}

\newcommand{\fraca}[2]{\displaystyle\frac{#1}{#2}}

\newcommand{\mm}[3]{\renewcommand{\arraystretch}{0.8}\begin{array}[t]{c}\mbox{#1}
\\ #2\end{array}\begin{array}[t]{c}#3\end{array}
\renewcommand{\arraystretch}{1}}

\def \R {{\rm I\kern -2.2pt R\hskip 1pt}}

\newtheorem{algorithm1}{Algorithm}[section]

\newtheorem{theorem}{Theorem}[section]
\newtheorem{proof}{Proof}[section]

\newtheorem{remark}{Remark}[section]

\def\boxforqed{\rule{0.5em}{1.5ex}}
\def\qed{\ifmmode\squareforqed\else{\unskip\nobreak\hfil
        \penalty50\hskip1em\null\nobreak\hfil\boxforqed
         \parfillskip=0pt\finalhyphendemerits=0\endgraf}\fi}

\begin{document}

\begin{frontmatter}

% Title, authors and addresses

% use the thanksref command within \title, \author or \address for footnotes;
% use the corauthref command within \author for corresponding author footnotes;
% use the ead command for the email address,
% and the form \ead[url] for the home page:
% \title{Title\thanksref{label1}}
% \thanks[label1]{}
% \author{Name\corauthref{cor1}\thanksref{label2}}
% \ead{email address}
% \ead[url]{home page}
% \thanks[label2]{}
% \corauth[cor1]{}
% \address{Address\thanksref{label3}}
% \thanks[label3]{}

\title{Robust solutions of uncertain mixed-integer linear programs using decomposition techniques}

% use optional labels to link authors explicitly to addresses:
% \author[label1,label2]{}
% \address[label1]{}
% \address[label2]{}

\author[label1]{R.~M\'{\i}nguez\corref{cor1}}
\address[label1]{Dr. Eng, HIDRALAB INGENIER\'IA Y DESARROLLOS, S.L., Spin-Off UCLM, Hydraulics Laboratory Univ.\ of Castilla-La Mancha, Av. Pedriza, Camino Moledores s/n, 13071 Ciudad Real, Spain}
\author[label2]{V.~Casero-Alonso}
\address[label2]{Department of Mathematics, Institute of Mathematics Applied to Science and Engineering, University of Castilla-La Mancha,
Ciudad Real, Spain}
\cortext[cor1]{Corresponding author: roberto.minguez@hidralab.com}

\begin{abstract}
% Text of abstract
Robust optimization is a framework for modeling optimization problems involving data uncertainty and during the last decades has been an area of active research. If we focus on linear programming (LP) problems with i) uncertain data, ii) binary decisions and iii) hard constraints within an ellipsoidal uncertainty set, this paper provides a different interpretation of their robust counterpart (RC) inspired from decomposition techniques. This new interpretation allows the proposal of an ad-hoc decomposition technique to solve the RC problem with the following advantages: i) it improves tractability, specially for large-scale problems, and ii) it provides the exact probability of constraint violation in case the probability distribution of uncertain parameters are completely defined by using first and second-order probability moments. An attractive aspect of our method is that it decomposes the second-order cone programming problem, associated with the robust counterpart, into a linear master problem and different quadratically constrained problems (QCP) of considerable lower size. The optimal solution is achieved through the solution of these master and subproblems within an iterative scheme based on cutting plane approximations of the second-order cone constraints. In addition, proof of convergence of the iterative method is given.
%
%A comprehensive case study based on the self-scheduling problem of
%a power producer is used to illustrate the capabilities of the
%proposed  methodology. Appropriate conclusions are finally drawn.
\end{abstract}

\begin{keyword}
% keywords here, in the form: keyword \sep keyword
Stochastic programming \sep Conic programming
and interior point methods \sep Decision analysis under uncertainty \sep Reliability analysis \sep Robust optimization
% PACS codes here, in the form: \PACS code \sep code

%% MSC codes here, in the form: \MSC code \sep code
%% or \MSC[2008] code \sep code (2000 is the default)

\end{keyword}
\end{frontmatter}

% main text
\section{Introduction}
%
%\subsection{Motivation}
%
The concept of robust optimization was developed to drop the classical assumption in mathematical programming that the input data is precisely known and equal to given nominal values. It is well known that in practice, most of the data involved in optimization problems is uncertain, and optimal solutions using nominal values might no longer be optimal or even feasible. Robust optimization techniques deal with the problem of designing solutions that are immune to data uncertainty \citep{Soyster:73,El-GhaouiL:97,El-GhaouiOL:98,Ben-TalN:98,Ben-TalN:99,Ben-TalN:00,BertsimasS:04} by solving equivalent deterministic problems. The main advantage of these techniques is that it is not required to know the probability density function (PDF) of the uncertain data. The decision-maker searches for the optimal solution of all-possible realizations of uncertain data within the uncertainty set, and in addition, probabilistic bounds of constraint violation valid for different probability density functions are available.

Stochastic programming is also a framework for modelling problems that involve uncertainty \citep{BirgeL:97}. In this particular case, uncertain data is assumed to follow a given probability distribution and are usually dealt with by using scenario models or finite sampling from the PDFs \citep{RockafellarU:00,RockafellarU:02}. However, the number of scenarios needed to represent the actual stochastic processes can be very large, which may result in intractable problems. That is the most important reason why robust optimization is gaining popularity among practitioners with respect to stochastic programming \cite{GabrelMT:14}, not only in the operational research community but also for design engineers \cite{HoushOS:11,PerelmanHO:13}.

Stochastic programming in the context of engineering design and optimization, i.e. reliability-based structural optimization \citep{Frangopol:95,Melchers:99,RoysetDP:01,RoysetDP:01b,RoysetDP:06}, has also broadened using as risk measure the failure probability. In this context, it is required to know: i) the joint probability density function of all random variables involved and ii) a method for calculating the probabilities of failure for a given design. Since the evaluation of failure probabilities is not an easy task, different methods have been developed, such as First-Order Second Moment (FOSM, \cite{HasoferL:74}).

Despite the analogies among the problems treated within these different frameworks, i.e. stochastic programming, robust optimization, and reliability-based structural optimization, research trends and solution techniques have followed different paths. To the best of our knowledge, a few works have taken advantage of methods from one field to be applied to any other. For instance, the work \citep{MinguezCG:11} proposes a new method to solve certain classes of stochastic programming problems based on FOSM and mathematical programming decomposition techniques. Their method focus on a specific type of problems where: i) the joint probability distribution of the random variables involved is given or can be estimated parametrically, ii) distributions do not depend on the decision variables, and iii) the random variables only affect the objective function. Recently, in \cite{RockafellarR:10} it is proposed a new risk measure, the buffered failure probability, which allows the generalization of the CVaR concept from stochastic programming \cite{RockafellarU:00,RockafellarU:02} to reliability-based structural optimization problems using finite sampling. One of the aims of the present work is to give a new perspective and apply concepts originating from structural reliability to robust optimization, we attempt to shed new light on existing problems and as such stir innovative thinking.

In particular, we focus on the type of problems dealt with on work \citep{Ben-TalN:99}, i.e. linear mathematical programming problems with hard constraints that must be satisfied for any possible realization of the uncertain data.  In paper \citep{Ben-TalN:99}, authors propose to obtain robust solutions of an uncertain LP problem with ellipsoidal uncertainty sets, whose RC results is a conic quadratic problem, i.e. a convex and tractable problem that can be solved in polynomial time by interior point algorithms. However, the inclusion of binary and/or integer decisions poses new challenges from the computational perspective.
%%%, which resulted in the proposal of an alternative uncertainty set that transforms the RC into a linear formulation solvable using standard branch-and-cut solvers \cite{BertsimasS:04}.
This paper proposes an alternative and decomposable solution technique based on cutting planes that allows reaching the optimal solution of the RC problem by solving two kind of problems within an iterative scheme: one mixed-integer linear master problem, and one subproblem of considerable lower size for each hard constraint. This strategy of decomposing a problem into smaller pieces has proved to be effective to improve tractability in many different applications \citep{Floudas:95,ConejoCMG:06}. Note that cutting plane algorithms for robust mixed-integer linear programs are state-of-the-art, see \citep{FischettiM:12,BertsimasDL:16} among many others. In particular, paper \citep{BertsimasDL:16} proposed a similar decomposition algorithm with respect to the one presented in this work with different variants to improve computational efficiency, however, proof of convergence and probabilistic guarantees are not given. Our research was conducted independently of the work in \citep{BertsimasDL:16}. There are also cutting-plane approximations with application to chance-constrained problems \citep{AckooijFD:16}.

The proposed method has the following features which makes it attractive for practical use:
%%%it improves tractability with respect to interior point algorithms when binary and/or integer decision variables are included, specially for large-scale problems, because
i) the master problem remains linear and ii) the subproblems are QCP with just one quadratic constraint, which have analytical solutions. In addition and due to the relationship between the subproblem formulation and reliability-based structural techniques (FOSM), it allows to calculate exact probabilities of constraint violation in case the probability distributions of uncertain parameters are completely defined by using first and second-order moments (mean and variance-covariance). This feature could encourage engineers used to working with failure probabilities to take advantage of robust optimization techniques, even without using the proposed iterative method. %Finally, and since in many applications it is not enough to know first and second moments of the probability %distributions, some hints about possible extensions of the method to work with current and more general multivariate %probabilistic models used by engineers are given.
%%%Finally, the method proposed in this paper constitutes an alternative formulation to that proposed by \cite{BertsimasS:04} using cardinality constrained uncertainty sets. Both approaches solve mixed-integer linear programming problems, the main difference is that the method proposed by \cite{BertsimasS:04} requires the definition of the maximum number of parameters allowed to reach its worst possible values, while the approach proposed in this paper requires the definition of the protection level. Besides, the proposed method would not require the definition of bounds for the uncertain parameters, which might be of interest for certain applications.

The rest of the paper is organized as follows. Section~\ref{s1} introduces the robust formulation of linear programming problems. In
Section~\ref{s.rborm} a detailed description of the decomposition method proposed in this work is given,
 while in Section~\ref{s2} the method for calculating probabilities of constraint violation using FOSM technique is explained in detail.
 %%%, paying special attention to its relationship with ellipsoidal uncertainty sets.
 %Discussion about the extension of the method to take into account different multivariate probability distribution %models is given in Section~\ref{}.
 In Sections~\ref{s.ilusexam} and \ref{s.sheexam} an illustrative example and a realistic case study are respectively
described, solved, analyzed, and compared with existing approaches. Section~\ref{discussion} briefly discusses future possible extensions of the method proposed in this work.
Finally, in Section~\ref{Conclu} some relevant conclusions are duly drawn.

\section{Robust counterpart of an uncertain mixed-integer linear programming problem}\label{s1}
Consider the following problem:
\begin{eqnarray}
\mm{Maximize}{{\bfg x}}{{\bfg c}^T{\bfg x}},&&\label{eq.fobj}\\
% \nonumber to remove numbering (before each equation)
 \mbox{subject to \hspace*{0.5cm}}{\bfg A} {\bfg x} &\le & {\bf 0} \label{eq.incons}\\
{\bfg l}\le {\bfg x} & \le &{\bfg u}\label{eq.bounds}\\
x_j &\in &{\cal Z};\; j=1, \ldots, k;\;k>0;\;k\le n,\label{eq.integrality}
\end{eqnarray}
where ${\bfg x}$ is the decision variable vector of $n$ variables, the first $k$ of which are integral. Note that we consider that $k>0$, i.e. it always contains binary and/or integer decisions, otherwise the problem can be solved efficiently using interior point algorithms. ${\bfg c}(n \times 1)$ and
${\bfg A}(m \times n)$ are data coefficients, and ${\bfg l}(n \times 1)$ and ${\bfg u}(n \times 1)$ are lower and upper decision variable bounds. We assume without loss of generality that the only uncertain coefficients are those belonging to matrix ${\bfg A}(m \times n)$. For those cases where vector ${\bfg c}(n \times 1)$ is uncertain, or even the right hand side of equation~(\ref{eq.incons}) is uncertain and equal to ${\bfg b}(m \times 1)$, it is possible to rewrite the original problem as (\ref{eq.fobj})-(\ref{eq.integrality}) (see \cite{Ben-TalEN:09}).

The RC of problem (\ref{eq.fobj})-(\ref{eq.integrality}) is the same problem but replacing constraint set (\ref{eq.incons}) by:
\begin{equation}\label{eq.inconsRC}
     {\bfg A} {\bfg x}  \le  {\bf 0};\;\forall  {\bfg A}\in {\bfg U},
\end{equation}
where ${\bfg U}$ is an uncertainty set. According to \cite{Ben-TalEN:09} (check also reference \cite{HoushOS:11}), an LP with a certain objective is a constraint-wise problem and its solution does not change if the uncertainty set is extended to the product of its projections on the subspaces of the constrains, i.e. constraint set (\ref{eq.inconsRC}) is equivalent to:
\begin{equation}\label{eq.inconsRCi}
     {\bfg a}^{(i)^T} {\bfg x}  \le  0;\;\forall  {\bfg a}^{(i)}\in U_i, i=1,\ldots,m,
\end{equation}
where ${\bfg a}^{(i)^T}; i=1,\ldots,m$ corresponds to the rows of matrix ${\bfg A}$ and $U_i$ is the projection of ${\bfg U}$ on the subspace of the data of ${\bfg a}^{(i)}$.

Traditionally, parameter uncertainty $a_{ij}\in {\bfg A}$ within robust optimization is modeled as a symmetric and bounded random variable $\tilde a_{ij}$ that takes values in the interval $[a_{ij}-\delta_{ij},a_{ij}+\delta_{ij}]$ following an unknown probability distribution. Elements $a_{ij}\in {\bfg A}$ represent nominal values and
$\delta_{ij};i=1,\ldots,m;j=1,\ldots,n$ are the maximum absolute value deviations from the nominal values. This interval characterization of parameter uncertainty is required if worst-case oriented methods, or box uncertainty sets are used. This is the solution proposed by \cite{Soyster:73}, where each uncertain parameter $\tilde a_{ij};i=1,\ldots,m;j=1,\ldots,n$ takes its worst possible value within the given interval. This strategy leads to an excessively conservative solution.

To address this excessive conservatism, paper \cite{Ben-TalN:99} proposes an alternative uncertainty set. Let us assume that uncertain parameter vectors $\tilde{\bfg a}^{(i)};i=1,\ldots,m$ have nominal or expected values ${\bfg a}^{(i)};i=1,\ldots,m$ and variance-covariance matrix ${\bfg \Sigma}^{(i)};i=1,\ldots,m$ which is positive definite, respectively. According to \cite{Ben-TalN:99}, the ellipsoidal uncertainty set can be written using the Mahalanobis distance as follows:
\begin{equation}\label{eq.Mahalanobis}
     U_i(\beta_i)\equiv\left\{\tilde{\bfg a}^{(i)}|\left(\tilde{\bfg a}^{(i)}-{\bfg a}^{(i)}\right)^T \left({\bfg \Sigma}^{(i)}\right)^{-1}\left(\tilde{\bfg a}^{(i)}-{\bfg a}^{(i)}\right)\le \beta_i^2 \right\};\; i=1,\ldots,m,
\end{equation}
so that the RC of problem (\ref{eq.fobj})-(\ref{eq.integrality}) is the same problem but replacing constraint set (\ref{eq.incons}) by:
\begin{equation}\label{eq.inconsRCellip}
   \left(\mm{Maximum}{\tilde{\bfg a}^{(i)}\in U_i(\beta_i)}{\tilde{\bfg a}^{(i)^T} {\bfg x}}\right) \le  0,\; i=1,\ldots,m,
\end{equation}
where parameters $\beta_i;\;i=1,\ldots,m$ control the size and protection level of the ellipsoidal sets for each constraint.

%%%According to \cite{Ben-TalN:99}, for given values of the decision variables ${\bfg x}$, the probability that each uncertain constraint is violated is given by the following expression:
%%%%
%%%\begin{equation}\label{eq.probound}
%%%   \mbox{Prob}\left(\tilde{\bfg a}^{(i)^T} {\bfg x}>0\right) \le  e^{-\beta_i^2/2}, i=1,\ldots,m.
%%%\end{equation}

Note that contrary to the worst case approach proposed by \cite{Soyster:73}, only first and second order moments of the random parameters without lower and upper bounds are considered in this paper for the ellipsoidal uncertainty set. The reason is that in case the uncertain parameters follow a multivariate normal distribution, the exact probability of constraint violation can be calculated, which might be of interest for practitioners.
Robust optimization approaches involving first and/or second-order moment information has been
studied by quite a few researchers (see for example \citep{DelageY:10}, which deals with distributionally robust optimization, and the recent comprehensive survey \citep{BertsimasBC:11} for
more references).
%%%, and secondly, it is straightforward to include also the interval limitations $[a_{ij}-\delta_{ij},a_{ij}+\delta_{ij}]$ of the random parameters within the proposed solution method, although tighter bounds with respect to those in (\ref{eq.probound}) would required a more sophisticated development. We will come back to these matters in the next sections.

In order to solve constraints (\ref{eq.inconsRCellip}), the uncertainty sets (\ref{eq.Mahalanobis}) are transformed using an affine mapping into balls of radius $\beta_i;i=1,\ldots,m$, respectively, resulting in the following set of alternative constraints:
\begin{equation}\label{eq.inconsRCellip2}
   \left(\mm{Maximum}{{\bfg z}^{(i)}}{\left({\bfg a}^{(i)}+ {\bfg L}^{(i)}{\bfg z}^{(i)} \right)^T {\bfg x}}\right) \le  0;\;\|{\bfg z}^{(i)}\|\le \beta_i, i=1,\ldots,m,
\end{equation}
where ${\bfg L}^{(i)}$ is the mapping matrix which can be obtained from Cholesky decomposition of variance-covariance matrix ${\bfg \Sigma}^{(i)}={\bfg L}^{(i)}{\bfg L}^{(i)^T}$,
${\bfg z}^{(i)}$ represents a perturbation vector and $\|\cdot\|$ stands for Euclidean norm. The analytical solution of constraints (\ref{eq.inconsRCellip2}) in terms of the objective function is (see \cite{Ben-TalN:99}):
\begin{equation}\label{eq.inconsRCellip3}
   {\bfg a}^{(i)^T}{\bfg x}+ \beta_i\|{\bfg L}^{(i)^T}{\bfg x} \| \le  0; i=1,\ldots,m,
\end{equation}
or equivalently,
\begin{equation}\label{eq.inconsRCellip3b}
  {\bfg a}^{(i)^T}{\bfg x}+\beta_i\sqrt{{\bfg x}^T{\bfg \Sigma}^{(i)}{\bfg x}} \le  0; i=1,\ldots,m.
\end{equation}

Finally, the RC of problem (\ref{eq.fobj})-(\ref{eq.integrality}) using ellipsoidal uncertainty sets becomes:
\begin{eqnarray}
&&\mm{Maximize}{{\bfg x}}{{\bfg c}^T{\bfg x}},\label{eq.fobjF}\\
% \nonumber to remove numbering (before each equation)
  %\sum_{j=1}^n a_{ij}x_j+\beta_i\sqrt{\sum_{k=1}^n\left(\sum_{j=k}^n L^{(i)}_{jk}x_j\right)^2} & \le & 0;\;i=1,\ldots,m,\label{eq.inconsF}\\
&&\mbox{\hspace*{0.2cm}subject to      }  {\bfg a}^{(i)^T}{\bfg x}+\beta_i\sqrt{{\bfg x}^T{\bfg \Sigma}^{(i)}{\bfg x}}  \le   0; i=1,\ldots,m,\label{eq.inconsF}\\
 % l_j\le x_j & \le & u_j;\;j=1,\ldots,n\label{eq.boundsF}\\
&& \hspace*{2.2cm} {\bfg l}\le {\bfg x}  \le  {\bfg u}\label{eq.boundsF}\\
&&  \hspace*{2.2cm} x_j  \in  {\cal Z};\; \forall j\le k;\;k>0;\;k\le n.\label{eq.integralityF}
\end{eqnarray}

Problem (\ref{eq.fobjF})-(\ref{eq.integralityF}) is a mixed-integer conic quadratic problem, i.e a nonlinear formulation, and as pointed out by \citep{BertsimasS:04}, it is not particularly attractive for solving robust linear discrete optimization models. Note that \citep{BertsimasS:04} proposed an alternative linear RC problem to avoid (\ref{eq.fobjF})-(\ref{eq.integralityF}) formulation when binary and/or integer variables are involved ($k>0$), which instead of allowing all random parameters to take their worst possible values within the given intervals such as \cite{Soyster:73}, only allows a pre-established number of parameters $\Gamma$ to reach them (polyhedral uncertainty sets). This alternative formulation remains linear, and it also provides a robust solution in terms of probability of infeasibility.

%%%%%%%%%%%%%%%%%%%%%%%%%%%%%%%%%%%%%%%%%%%%%%%%%%%%%%%%%%%%%%%%%%%%%%%%%%%%%%%%%%%%%%%%%%%%%%%%%%%%%%%%%%%%%%%%
%%%%%%%%%%%%%%%%%%%%%%%%%%%%%%%%%%%%%%%%%%%%%%%%%%%%%%%%%%%%%%%%%%%%%%%%%%%%%%%%%%%%%%%%%%%%%%%%%%%%%%%%%%%%%%%%
%%%%%%%%%%%%  INICIO DE LA PARTE DE VICTOR
%%%%%%%%%%%%%%%%%%%%%%%%%%%%%%%%%%%%%%%%%%%%%%%%%%%%%%%%%%%%%%%%%%%%%%%%%%%%%%%%%%%%%%%%%%%%%%%%%%%%%%%%%%%%%%%%
%%%%%%%%%%%%%%%%%%%%%%%%%%%%%%%%%%%%%%%%%%%%%%%%%%%%%%%%%%%%%%%%%%%%%%%%%%%%%%%%%%%%%%%%%%%%%%%%%%%%%%%%%%%%%%%%

\section{Decomposition method for the RC}\label{s.rborm}
This paper provides a methodology to solve problem (\ref{eq.fobjF})-(\ref{eq.integralityF}) using mathematical programming decomposition techniques.
%Let assume ${\bfg x}_\nu;\nu=1,\ldots,\infty$ being different values of the design variables ${\bfg x}$, thus
The mentioned problem can be rewritten equivalently as follows:
%
%\begin{equation}\label{eq.fobjF2}
   % \mm{Maximize}{x_j;j=1,\ldots,n}{\displaystyle\sum_{j=1}^n c_jx_j},
%\end{equation}
%
%subject to
%
\begin{eqnarray}
% \nonumber to remove numbering (before each equation)
&&\mm{Maximize}{{\bfg x}}{{\bfg c}^T{\bfg x}},\label{eq.fobjF2}\\
&& \hspace*{0.2cm}\mbox{subject to } \hat {\bfg a}_{\nu}^{(i)^T}{\bfg x} \le 0;\;i=1,\ldots,m;\;\nu=0,\ldots,\infty, \label{eq.inconsF2}\\
  %l_j\le x_j & \le & u_j;\;j=1,\ldots,n,\label{eq.boundsF2}\\
 && \hspace*{2.2cm}{\bfg l}\le {\bfg x} \le  {\bfg u}\label{eq.boundsF2}\\
 && \hspace*{2.2cm}x_j  \in  {\cal Z};\; \forall j\le k;\;k>0;\;k\le n,
  \label{eq.integralityF2}
\end{eqnarray}
%
%%%\footnote{$\hat {\bfg a}_\nu^{(i)^T}{\bfg x}$  or $\sum_{j=1}^n \hat a_{ij_\nu} x_j$}
where column-vectors $\hat {\bfg a}_{\nu}^{(i)}=(\hat a_{i1_{\nu}}, \ldots,\hat a_{ij_{\nu}},\ldots,\hat a_{in_{\nu}})^T;i=1,\ldots,m$ for constraint $\nu$ correspond, respectively, to the solution of the following optimization problems:
\begin{equation}\label{FOSMEquivalent2}
\left.
\begin{array}{l}
\hat {\bfg a}^{(i)}_\nu = \mm{arg max}{ \tilde {\bfg a}^{(i)}}{\tilde{\bfg a}^{(i)^T} {\bfg x}_{\nu}},\\
\mbox{subject to}\\
%%% \tilde {\bfg a}^{(i)} &= &{\bfg a}^{(i)} + {\bfg L}^{(i)}{\bfg z}^{(i)},\\
%%%\sqrt{\sum\limits_{j=1}^{n}z_{j}^{(i)^2}}& \le &\beta_i,
\left(\tilde{\bfg a}^{(i)}-{\bfg a}^{(i)}\right)^T \left({\bfg \Sigma}^{(i)}\right)^{-1}\left(\tilde{\bfg a}^{(i)}-{\bfg a}^{(i)}\right)\le \beta_i^2,
\end{array}
\right\}
%\begin{array}{l}
%i=1,\ldots,m,\\
\nu=1,\ldots,\infty,
%\end{array}
\end{equation}
where ${\bfg x}_{\nu}$ is the solution of (\ref{eq.fobjF2})-(\ref{eq.integralityF2}) when constraints $0,\ldots,\nu-1$ are considered, and $\hat {\bfg a}^{(i)}_0$ are the expected or nominal values of the problem, ${\bfg a}^{(i)};\forall i$.

Problem (\ref{FOSMEquivalent2}) is equivalent to the following problem:
\begin{equation}\label{FOSMEquivalent3}
\left.
\begin{array}{l}
\hat {\bfg a}^{(i)}_\nu = \mm{arg max}{ \tilde {\bfg a}^{(i)}}{\tilde{\bfg a}^{(i)^T} {\bfg x}_{\nu}},\\
\mbox{subject to}\\
 \tilde {\bfg a}^{(i)} ={\bfg a}^{(i)} + {\bfg L}^{(i)}{\bfg z}^{(i)},\\
 {\bfg z}^{(i)^T}{\bfg z}^{(i)}\le \beta_i
\end{array}
\right\}i=1,\ldots,m,
\end{equation}
which corresponds to the problem defined in constraint (\ref{eq.inconsRCellip2}) particularized for decision variables ${\bfg x}_{\nu}$.

In this alternative %the so-called points of maximum likelihood for each constraint
the values of $\hat {\bfg a}_\nu^{(i)}$ %(we use this name from the relationship with FOSM methods, as explained in next section)
are used explicitly in (\ref{eq.inconsF2}) to define tangent hyperplanes with respect to the original convex conic restrictions.
%%%This is because the structure of the problem (\ref{FOSMEquivalent2}), whose solution is precisely the derivative of the conical restriction
%%%(as we will see in Remark \ref{analyticalsol.remark}).

\begin{remark}
Theoretically, the conical convex constraints (\ref{eq.inconsF}) can be reproduced with an infinite number of tangent hyperplanes defining a linear envelope, as described in (\ref{eq.inconsF2}).

%This makes problems (\ref{eq.fobjF})-(\ref{eq.integralityF}) and (\ref{eq.fobjF2})-(\ref{FOSMEquivalent2})
%to be equivalents, as shown in Theorem~\ref{Equiv.th}
%However, the iterative process above converges to the optimal solution by using a finite number of hyperplanes (Theorem \ref{Equiv.th}).
\end{remark}

%Note that constraint (\ref{eq.inconsF2}) constitute a linear envelope of the original convex constraints (\ref{eq.inconsF}) composed theoretically by infinite hyperplanes, which makes problems (\ref{eq.fobjF})-(\ref{eq.integralityF}) and (\ref{eq.fobjF2})-(\ref{FOSMEquivalent2}) to be equivalents, as shown in Theorem~\ref{Equiv.th}. However,

The question is: how to define efficiently those possible solutions ${\bfg x}_\nu$ in order to construct cuts (\ref{eq.inconsF2}) without the need to include an unlimited number of constraints? To answer this question, we propose decomposing the original problem into two procedures:
\begin{enumerate}
  \item {\bf Decision making at iteration $l$:} For given values of %the points of maximum likelihood
  $\hat {\bfg a}^{(i)}_\nu;$ $i=1,\ldots,m;$ $\nu=0,...,l-1,$ the decision variables maximizing the problem (\ref{eq.fobjF2})-(\ref{eq.integralityF2})  are obtained, i.e.  $\bfg x_l$. This is considered the {\it master problem}. Due to the convex character of conic constraints, the master problem is a relaxation of the original problem because the feasible region defined by the linear envelope always contain the feasible region defined by the original conic restrictions.
  % (as we will see in Remark \ref{analyticalsol.remark}).
   Note that the master problem solution constitutes an upper bound of the optimal solution.
  \item {\bf Construct additional tangent hyperplanes:} For the decisions $\bfg x_l$ made in the previous step, update the values of the random variables $\hat {\bfg a}^{(i)}_l$ required to achieve the target security criterium given by ${\bfg \beta}=(\beta_1,\beta_2,\ldots,\beta_m)^T$, which is equivalent to construct additional tangent hyperplanes with respect to the original conic constraints. These are considered the {\it subproblems}.
\end{enumerate}

Thus, the solution of the alternative problem (\ref{eq.fobjF2})-(\ref{FOSMEquivalent2}) is achieved by means of an iterative scheme, which is repeated until the stopping criterion is satisfied.
%%%or \textst{the values of the decision variables $\bfg{x}$ stabilize} almost, i.e., when a tolerance $\varepsilon$ instead the zero is considered in (\ref{eq.inconsF2}).
The main reasons for proposing this decomposition are:
\begin{enumerate}
  \item Problem (\ref{eq.fobjF2})-(\ref{eq.integralityF2}) for given values of the parameter values $\hat {\bfg a}^{(i)}_\nu$ is a mixed-integer linear mathematical programming problem, which can be solved efficiently using state-of-the-art mixed-integer solvers.

  \item Problems (\ref{FOSMEquivalent2}), for given values of the decision variables ${\bfg x}_\nu$, can be solved for each constraint independently. Each of these problems is a QCP problem with just one quadratic constraint \citep{Pinar:00}, which has an analytical solution as shown below.
\end{enumerate}

%\begin{lemma}[Analytical solution of subproblem (\ref{FOSMEquivalent2})]\label{analyticalsol}
%Given any solution $x_\nu$, an hyperplane can be construct solving the subproblem (\ref{FOSMEquivalent2})
The analytical optimal solutions of subproblems (\ref{FOSMEquivalent3}) at iteration $l$ for given values of the decision variables ${\bfg x}_l$, in terms of the uncertain coefficients, are:
 \begin{equation}\label{eq.ahati}
    \hat {\bfg a}^{(i)}_l={\bfg a}^{(i)}+\beta_i\fraca{{\bfg \Sigma}^{(i)}{\bfg x}_l}{\sqrt{{\bfg x}_l^T{\bfg \Sigma}^{(i)}{\bfg x}_l}};\;i=1,\ldots,m.
\end{equation}

This result is straightforward to verify since the optimal solution of convex problem (\ref{FOSMEquivalent3}) in terms of the objective function is given by (\ref{eq.inconsRCellip3b}) particularized for ${\bfg x}_l$,  which is precisely the optimal objective function obtained if the objective function in (\ref{FOSMEquivalent3}) is evaluated at the optimum (\ref{eq.ahati}). Nevertheless, a formal proof of this result is given in \cite[Lemma 1]{Ackooij:15}.

%\end{lemma}

%\begin{proof}[Analytical solution of subproblem (\ref{FOSMEquivalent2})]
%It is straightforward to see that (\ref{eq.ahati}) is the optimal solution of problem (\ref{FOSMEquivalent2}) by just evaluating its objective function at this point, i.e.:
%

%
%which is exactly the same optimal objective function as that given by expression (\ref{eq.inconsRCellip3b}).
%This completes the proof.
%\end{proof}

\begin{remark}\label{analyticalsol.remark}
Note that parameters $\hat {\bfg a}_l^{(i)}$ in (\ref{eq.ahati}) are the derivatives of the conic constraints (\ref{eq.inconsF}) at $\bfg{x}_l$:
\begin{equation}\label{eq.ahatiproof}
 \begin{array}{rcl}
   \hat {\bfg a}_l^{(i)^T}{\bfg x}_l & = & {\bfg a}^{(i)^T}{\bfg x}_l+\beta_i\fraca{{\bfg x}_l^T{\bfg \Sigma}^{(i)^T}{\bfg x}_l}{\sqrt{{\bfg x}_l^T{\bfg \Sigma}^{(i)}{\bfg x}_l}};\;i=1,\ldots,m,\\
    & = & {\bfg a}^{(i)^T}{\bfg x}_l+\beta_i\sqrt{{\bfg x}_l^T{\bfg \Sigma}^{(i)}{\bfg x}_l};\;i=1,\ldots,m,
 \end{array}
\end{equation}
which explains why constraints (\ref{eq.inconsF2}) correspond to tangent hyperplanes with respect to the conic restrictions. Then, the feasible region defined by the linear envelope always contain the feasible region defined by the original conic restrictions.
The latter can be shown as follows. For all feasible solution, $\bfg x_f$, of the original problem (\ref{eq.inconsF})-(\ref{eq.integralityF}) it is satisfied
$$  {\bfg a}^{(i)^T}{\bfg x}_f+\beta_i\sqrt{{\bfg x}_f^T{\bfg \Sigma}^{(i)}{\bfg x}_f} \leq  0;\;i=1,\ldots,m.$$
Then
$$ \hat {\bfg a}_{l}^{(i)^T} \bfg x_f = {\bfg a}^{(i)^T}{\bfg x}_f+\beta_i \fraca{{\bfg x}_l^T{\bfg \Sigma}^{(i)^T}{\bfg x}_f}{\sqrt{{\bfg x}_l^T{\bfg \Sigma}^{(i)}{\bfg x}_l}}
 \leq {\bfg a}^{(i)^T}{\bfg x}_f+\beta_i\sqrt{{\bfg x}_f^T{\bfg \Sigma}^{(i)}{\bfg x}_f} \le 0;\;i=1,\ldots,m$$
 for all $l$
 because ${\bfg x}_l^T{\bfg \Sigma}^{(i)^T}{\bfg x}_f \le
 \sqrt{{\bfg x}_f^T{\bfg \Sigma}^{(i)}{\bfg x}_f}
 \sqrt{{\bfg x}_l^T{\bfg \Sigma}^{(i)}{\bfg x}_l}
 $
due to Cauchy-Schwarz inequality.
Therefore, $\bfg x_f$ is a feasible solution of (\ref{eq.inconsF2})-(\ref{eq.integralityF2}).
\end{remark}

%%%
%%%\footnote{Quitamos entonces todo el parrafo?} We examine the ?size of formulations? \footnote{7.- Igual hay que replantearse lo de escribir la comparativa de número de restricciones y variables porque ahora no es fácil saber cuantas restricciones se van metiendo porque se acumulan.}to solve problem (\ref{eq.fobj})-(\ref{eq.integrality}) (initial problem) using the RC formulation (\ref{eq.fobjF})-(\ref{eq.integralityF}) proposed by \cite{Ben-TalN:99}, and the master (\ref{eq.fobjF2})-(\ref{eq.integralityF2}) and subproblems (\ref{FOSMEquivalent2}) proposed in this paper. Assuming that all nominal values in matrix ${\bfg A}$ are subject to uncertainty, and given that the original problem has $n$ variables and $m$ constraints (not counting decision variable bounds), (\ref{eq.fobjF})-(\ref{eq.integralityF}) is a mixed-integer second-order cone mathematical programming problem (QCP) with $n$ variables and $m$ quadratic constrains, the master problem (\ref{eq.fobjF2})-(\ref{eq.integralityF2}) is a mixed-integer linear mathematical programming problem with $n$ variables and $m$ constrains, and each of the $m$ subproblems is a QCP with $n$ variables, and one quadratic constraint.

The proposed iterative scheme to solve (\ref{eq.fobjF2})-(\ref{FOSMEquivalent2}) is described below:

\begin{algorithm1}{Decomposition method}\label{alg}

\begin{enumerate}
\item[Step 0] \textbf{Problem definition:} Selection of protection levels $\beta_i;\forall i$, objective function, constraints, the mean and variance-covariance matrix of the involved random parameters $\tilde{\bfg a}^{(i)};\forall i$ and the tolerance of the process $\varepsilon$. Set $l=0$ and the initial values of the random parameters to their expected or nominal values $\hat{\bfg a}^{(i)}_0= {\bfg a}^{(i)};\forall i$.
\item[Step 1] \textbf{Solving the master problem:} Update the iteration counter $l \longrightarrow l+1$ and calculate the optimal solution ${\bfg x}_l$ of the following master problem:
%\begin{equation}\label{eq.fobjF2.Step1}
%    \mm{Maximize}{{\bfg x}}{{\bfg c}^T{\bfg x}},
%\end{equation}
%
%subject to
%
\begin{eqnarray}
    \mm{Maximize}{{\bfg x}}{{\bfg c}^T{\bfg x}},\label{eq.fobjF2.Step1}\\
 \mbox{subject to \hspace*{0.1cm} }\hat {\bfg a}_\nu^{(i)^T}{\bfg x} & \le & 0;\;i=1,\ldots,m;\;\nu=0,1,\ldots,l-1 \label{eq.inconsF2.Step1}\\
  {\bfg l}\le {\bfg x} & \le & {\bfg u},\label{eq.boundsF2.Step1}\\
  x_j & \in & {\cal Z};\; \forall j\le k;\;k>0;\;k\le n.\label{eq.integralityF2.Step1}
\end{eqnarray}

Continue in {\bf Step 2}.

\item[Step 2] \textbf{Stopping rule:} Check if the current solution $\bfg{x}_l$ satisfies the original conic restrictions (\ref{eq.inconsF}). If it does, the optimal solution has been found. If it does not but $l>1$ and $\hat {\bfg a}_\nu^{(i)^T}{\bfg x}_l<\varepsilon$, stop the process with ${\bfg x}_l$ as optimal solution, otherwise continue to {\bf Step 3}.
\item[Step 3] \textbf{Solving subproblems:} Plug solution ${\bfg x}_l$ obtained in {\bf Step 1} in (\ref{eq.ahati}) to obtain $\hat{{\bfg a}}^{(i)}_l$ and continue in {\bf Step 1}.
\end{enumerate}
\end{algorithm1}
We iteratively repeat Steps from 1 to 3, until the stopping rule holds.
%%%optimal solution is achieved or the \textst{norm associated with design variables $||{\bfg x}_{l}||$ does not change significantly} the value of the tangent hyperplane at the last solution, $\bfg x_l$, is close to achieve the conical restriction (less than the tolerance $\varepsilon$) \textst{between consecutive iterations $l$ and $l-1$}.

\begin{remark}\label{remark33}
Theoretically, constraints (\ref{eq.inconsF2.Step1}) contain the hyperplanes associated with all rows of matrix ${\bfg A}$, i.e. $m$ tangent hiperplanes for each iteration. However, in practice, only the hyperplane approximations related to infeasible conic restrictions at the current solution point ${\bfg x}_l$ are required to achieve convergence, i.e. tangent hyperplanes related to $i\in \{1,2,\ldots,m\}$ that satisfies condition
$${\bfg a}^{(i)^T}{\bfg x}_l+\beta_i\sqrt{{\bfg x}_l^T{\bfg \Sigma}^{(i)}{\bfg x}_l} >  0,$$
are used.
\end{remark}

%\begin{remark} Moreover, through several examples considered in the development of this work, all with a optimal solution in a vertex, the solution of the proposed decomposition method is very efficient. In all these cases, only the last restriction (hyperplane) is needed to obtain a fast convergence to the optimal solution (case k=0 included).\end{remark}

 Let us remark the reader that the non-linearities in formulation (\ref{eq.fobjF})-(\ref{eq.integralityF}) do not make this problem particularly attractive  for solving robust discrete optimization models, however, the proposed decomposition remove those non-linearities from the master problem, which can be solved efficiently using standard mixed-integer linear programming algorithms.

The convergence characteristics of this iterative method are discussed in the next Theorem based on the ideas of the outer approximation algorithm (\citep{FletcherL:94, DuranG:86}) and the extended Kelly's cutting plane algorithm (\citep{WesterlundP:95}).
%%%\ref{apen.B}.
\begin{theorem} The Algorithm \ref{alg} using decomposition techniques
terminates in a finite number of steps
and solves the problem (\ref{eq.fobjF})-(\ref{eq.integralityF}).
\end{theorem}

\begin{proof}
%The proof is obtained directly from Theorem 2 of \citep{FletcherL:94}.
%Our decision space $X$ is a nonempty compact convex set defined by a system of linear inequality constraints. And the functions involved in our problem are convex and continuously differentiable.
%In our case, it is straightforward to see that assumptions A1 and A2 hold.
%Moreover \citep{Lobo:98} showed that the robust linear programming is convex.
%And \citep{Boyd:04} expressed the robust linear programming as a second-order cone constraint.
%The last condition before applying Theorem 2 of \citep{FletcherL:94} is the finiteness of the space of integer variables. This is true in our case because all the variables, continuous or discrete, are bounded.

We prove that no solution is generated twice by the iterative process. %follows the finiteness of the algorithm.
At iteration $l$, the solution $\bfg{x}_l$ for the master problem (\ref{eq.fobjF2.Step1})-(\ref{eq.integralityF2.Step1}) is obtained. If it is not the optimal solution, at least one of the original conic constraints (\ref{eq.inconsF}) is not satisfied, for instance ${\bfg a}^{(j)^T}{\bfg x}_l+\beta_j \sqrt{{\bfg x}_l^T{\bfg \Sigma}^{(j)}{\bfg x}_l} > 0$ for some $j\in\{1,\ldots,m\}$. Then, the corresponding tangent hyperplane $\hat{\bfg a}_l^{(j)^T} \bfg{x} \le 0$ is added for the next master problem, and since it is infeasible at iteration $l$ because $\bfg a^{(j)^T}{\bfg x}_l+\beta_j \sqrt{{\bfg x}_l^T{\bfg \Sigma}^{(j)}{\bfg x}_l} > 0$, the next solution of the master problem $\bfg{x}_{l+1}$ must be different so that $\hat{\bfg a}_l^{(j)^T} \bfg{x}_{l+1}\le 0$.

Then, the finiteness of Algorithm \ref{alg} follows from the previous property and from:
\begin{enumerate}
  \item The finiteness of the feasible region if all variables are integral.
  \item The convergence of the extended Kelley's cutting plane method for convex MINLP problems \citep{WesterlundP:95}. In reference \citep{WesterlundP:95} it is proved that for practical computation, replacing constraints $\hat {\bfg a}_\nu^{(i)^T}{\bfg x}  \le  0;\;i=1,\ldots,m;\;\nu=0,1,\ldots,l-1$ in (\ref{eq.fobjF2.Step1}) by $\hat {\bfg a}_\nu^{(i)^T}{\bfg x}  \le \varepsilon;\;i=1,\ldots,m;\;\nu=0,1,\ldots,l-1$, being $\varepsilon$ a given tolerance, ensures that convergence is clearly achieved in a finite number of steps. Note that the difference between the extended Kelley's cutting plane method \citep{WesterlundP:95} and the algorithm proposed in this paper consist of how the new tangent hyperplanes are defined. In our case we have an explicit formula, while the algorithm proposed by \citep{WesterlundP:95} requires checking the more restrictive convex restriction and computing its derivatives.
%In case continuous variables are present, which ends the iterative process when the norm associated with decision variables between two consecutive iterations changes less than the tolerance $\varepsilon$ or all conic constraints are feasible.
\end{enumerate}
Now it is shown that the proposed method always terminates at a solution of (\ref{eq.fobjF})-(\ref{eq.integralityF}). Let $\bfg{x}^\star$ be the optimal solution of (\ref{eq.fobjF})-(\ref{eq.integralityF}).
Since (\ref{eq.fobjF2.Step1})-(\ref{eq.integralityF2.Step1}) is a relaxation of (\ref{eq.fobjF})-(\ref{eq.integralityF}) (see Remark \ref{analyticalsol.remark}), ${\bfg c} ^T \bfg{x}^\star$ is a lower bound with respect to the optimal value of (\ref{eq.fobjF2.Step1})-(\ref{eq.integralityF2.Step1}), which is attained at $\bfg{x}^\star$. Now assume that $\bfg{x}_l$  is the solution of (\ref{eq.fobjF2.Step1})-(\ref{eq.integralityF2.Step1}) with ${\bfg c} ^T \bfg{x}_l<{\bfg c} ^T \bfg{x}^\star$ (i.e. not optimal), however, $\bfg{x}^\star$ must be feasible in the previous step, which contradicts the assumption that $\bfg{x}_l$ is the optimal solution of (\ref{eq.fobjF2.Step1})-(\ref{eq.integralityF2.Step1}). This concludes the proof.

\end{proof}
%%%\begin{remark}
%%%If all the variables considered are binary or integers, i.e. $x_i \in \mathcal{Z};\;i=1,\ldots,n$ or equivalently $k=n$, the finiteness of Algorithm \ref{alg} follows from the finiteness of those variables, that are bounded, and the property that no integer solution is generated twice by the algorithm. Note that assuming that the original problem is feasible, and if all the decision variables are binary/integer, the feasible region is a finite set (although it could be extremely large), so if our algorithm improves the solution between iterations and the number of solutions is finite, the algorithm would stop in a number of finite steps (which could also be extremely high).
%%%\end{remark}
\begin{remark}
We have not proved the convergence rate, however numerical simulations with different physically based problems indicate that convergence is achieved in a reduced number of iterations.
%%%In case the decision variable-vector contains continuous and discrete variables and if the problem is also feasible, we could fix the binary/integer variables to given values, so that the master problem solution is the result of solving a linear mathematical programming problem. If the problem is well posed and it has a unique solution in terms of continuous variables, we can thus assume that the possible solution space is conformed by a finite set of possible values of binary/integer variables, each of them with one optimal solution in terms of continuous variables which is always attained because branch-and-bound linear mathematical solvers are used for the master problem. Since our algorithm does provide a different solution between consecutive iterations, it would also attain convergence in a number of finite steps.
%Note that we use the same arguments than those used in the paper by Fletcher  and Leyffer \citep{FletcherL:94} to prove the convergence of solving mixed integer nonlinear programs by outer approximation, because our algorithm constitutes an outer approximation.
\end{remark}

Recently, paper \citep{BertsimasDL:16} proposed a similar algorithm to that proposed in this work. The difference with respect our method lies in the solution of our called master problem. We seek optimality in the solution at each iteration while this condition is relaxed in \citep{BertsimasDL:16} to gain computing speed. No proof of convergence of the algorithm is given in \citep{BertsimasDL:16}.

\section{Probability of constraint violation }\label{s2}
Once the optimal solution $\bfg{x}^\star$ of the robust problem (\ref{eq.fobjF})-(\ref{eq.integralityF}) is obtained, it might be of interest to calculate the probability of each constraint violation, i.e. $\mbox{Prob}({\tilde{\bfg a}}^{(i)^T}\bfg{x}^\star>0);i=1,\ldots,m$. To that end, we use the First-Order Second-Moment method from structural reliability \citep[see][]{Freudenthal:56,HasoferL:74,RackwitzF:78,Ditlevsen:81,HohenbichlerR:81,Melchers:99,MinguezC:09}, which requires
calculating the following parameter for each constraint:
\begin{equation}\label{FOSM}
\left.
\begin{array}{l}
 \hat \beta_i = \mm{Minimum}{{\bfg z}^{(i)}}{\sqrt{\sum\limits_{j=1}^{n}z_{j}^{(i)^2}}},\\
\mbox{subject to}\\
 \tilde {\bfg a}^{(i)} = {\bfg a}^{(i)} + {\bfg L}^{(i)}{\bfg z}^{(i)},\\
\tilde{\bfg a}^{(i)^T} {\bfg x}^\star=   0
\end{array}
\right\}i=1,\ldots,m,
\end{equation}
where the optimal solution $\hat{\bfg z}^{(i)}$ corresponds to the closest point to the origin located on the limit of constraint violation in the standard and independent normal random space, $\hat\beta_i$ is the minimum distance so-called {\em reliability index} in the structural reliability scientific community, and
$\hat{\bfg a}^{(i)}={\bfg a}^{(i)} + {\bfg L}^{(i)}\hat{\bfg z}^{(i)}$ is the {\em point of maximum likelihood}, i.e. the actual values of the uncertain parameters that make constraints to be active where the probability is higher, and it represents the most likely values of the random parameters that produce constraint violation.
 Note that the reliability index $\hat\beta_i$ is a non-negative value for probabilities of failure/infeasibility lower than $0.5$, which is the case for the robust application considered in this paper. The final probability of constraint violation is related to the reliability
index by the relation:
\begin{equation}\label{eq.pfdf}
   \mbox{Prob}\left({\tilde{\bfg a}}^{(i)^T}\bfg{x}^\star>0\right)=\Phi (-\hat\beta_i);\forall i=1,\ldots, m,
\end{equation}
where $\Phi (\cdot )$ is the cumulative distribution function of
the standard normal random variables. This method provides the
exact probability if the limit-state equation is linear in the standard
normal random space, i.e., if the resulting limit-state
distribution is normally distributed, which is the case if uncertain parameters are normally distributed.

From the practical perspective, the calculation of the reliability index $\hat\beta_i$ for each constraint from (\ref{FOSM}), needed to compute the probability of constraint violation, is only required for those inactive constraints at the optimal solution $\bfg{x}^\star$, otherwise its value corresponds to the selected protection level $\beta_i$.

According to (\ref{eq.pfdf}), it is worth pointing out that the robust constraints (\ref{eq.inconsRCi}) of an uncertain mixed-integer linear programming problem, assuming that random variables ${\tilde{\bfg a}}^{(i)};\;\forall i=1,\ldots,m$ follow multivariate gaussian distributions, can be
equivalently formulated as probabilistic or chance constraints as follows:
\begin{equation}\label{eq.inconsRCi2}
    \mbox{Prob}\left( {\tilde{\bfg a}}^{(i)^T} {\bfg x}  \le  0\right)\ge \Phi(\beta_i);\;\forall i=1,\ldots,m.
\end{equation}

\section{Illustrative example}\label{s.ilusexam}
In order to illustrate the proposed method and the graphical interpretation of the iterative process, a simple example  with only two decision variables is presented below.

Let consider the following problem:
\begin{eqnarray}
% \nonumber to remove numbering (before each equation)
&&\mm{Maximize}{x_1,x_2}{c_1x_1+c_2x_2},\label{eq.fobjIllusEx}\\
&&\hspace*{0.2cm}\mbox{subject to }a_{11}x_1+a_{12}x_2  \le  b_1,\label{eq.inconsF1IllusEx}\\
 && \hspace*{2.2cm} a_{21}x_1+a_{22}x_2  \le  b_2,\label{eq.inconsF2IllusEx}\\
 &&  \hspace*{2.2cm}   x_{i}  \in  \{0,1,2,\ldots,10\};\;i=1,2,\label{integralityIE}
\end{eqnarray}
where ${\bfg c}=(3\;3)^T$, ${\bfg A}=\left(\begin{array}{cc}1 & 2 \\2 & 1 \end{array}\right)$, and ${\bf b}=(7\;7)^T$. The only uncertain parameters are those of matrix ${\bfg A}$ so that their expected values are equal to their nominal values and the variance-covariance matrix associated with each constraint (\ref{eq.inconsF1IllusEx}) and (\ref{eq.inconsF2IllusEx}) are, respectively: $${\bfg \Sigma}^{(1)}={ \footnotesize \left(\begin{array}{cc}0.1^2 & 0.016 \\0.016 & 0.2^2                                                                 \end{array}\right)};\;{\bfg \Sigma}^{(2)}={ \footnotesize \left(\begin{array}{cc}0.2^2 & -0.01 \\-0.01 & 0.1^2                                                                 \end{array}\right).}$$

\begin{figure}[t]
\begin{center}
 \includegraphics*[width=0.95\textwidth]{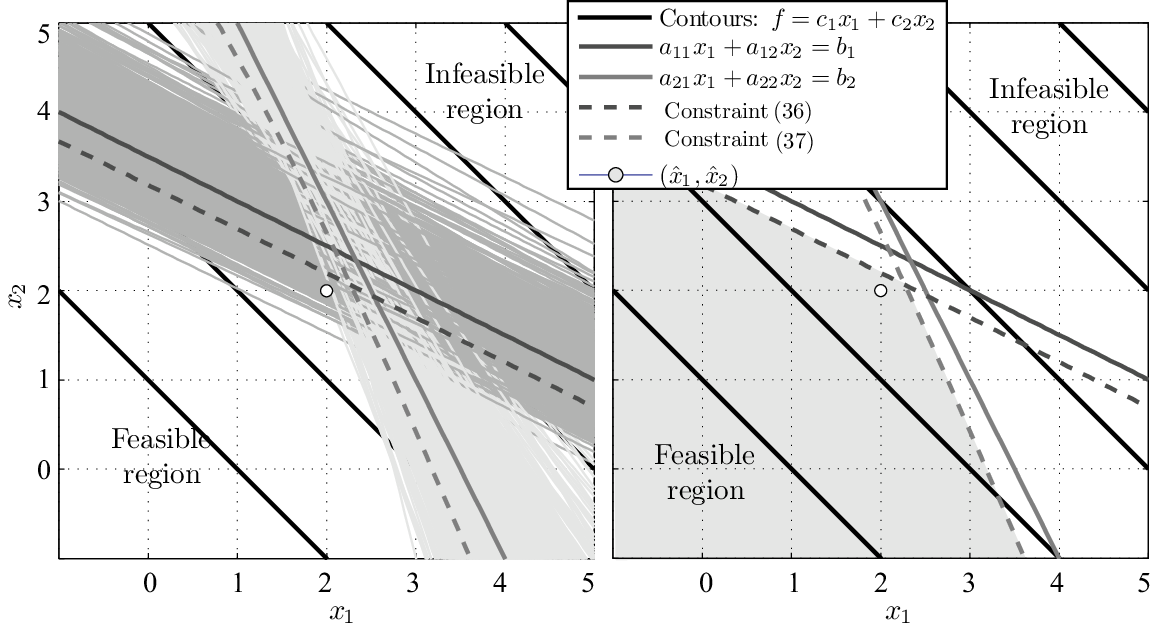}
\caption{\label{EjemploIlustrativo1Editado} Graphical illustration of illustrative example: panel on the left shows random instances of uncertain constraints, while panels on the right show the conic robust constraints that define the feasibility region and the optimal solution.}
\end{center}
\end{figure}

Figure~\ref{EjemploIlustrativo1Editado} (panel left) shows the graphical illustration of the problem (\ref{eq.fobjIllusEx})-(\ref{integralityIE}). The feasible region of the nominal problem is defined by the two gray constraints, while the black lines are contours associated with different values of the objective function. Note that its value is higher inside the unfeasible region. The gray and light gray shadows are indeed 1000 different realizations of the uncertain constraints, and it can be observed that given the optimal solution (white circle ($\hat x_1,\hat x_2$)), this point is unfeasible for many realizations of the uncertain constraints.

The RC of problem  (\ref{eq.fobjIllusEx})-(\ref{integralityIE}) using \citep{Ben-TalN:99} formulation is:
\begin{eqnarray}
% \nonumber to remove numbering (before each equation)
&& \mm{Maximize}{x_1;x_2}{c_1x_1+c_2x_2},\label{eq.fobjFIllusEx2}\\
&&\hspace*{0.2cm}\mbox{subject to }a_{11}x_1+a_{12}x_2 +\beta_1\sqrt{\sum_{k=1}^2\left(\sum_{j=k}^2 L^{(1)}_{jk}x_j\right)^2}  \le  b_1,\label{eq.inconsF1IllusEx2}\\
&&  \hspace*{2.2cm} a_{21}x_1+a_{22}x_2 +\beta_2\sqrt{\sum_{k=1}^2\left(\sum_{j=k}^2 L^{(2)}_{jk}x_j\right)^2} \le  b_2,\label{eq.inconsF2IllusEx2}\\
&&   \hspace*{2.2cm}      x_{i} \in  \{0,1,2,\ldots,10\};i=1,2,\label{integralityIE2}
\end{eqnarray}
which assuming $\beta_1=\beta_2=1$, results in the following optimal solution: $\hat f=10$, $\hat x_1= 2$, $\hat x_2= 2$. Once the optimal values of the decision variables are obtained we check the true values of the reliability indexes by solving problem (\ref{FOSM}). They are equal to $\hat \beta_1=1.746$ and $\hat \beta_2=2.887$ confirming that they are above the minimum required values of $1.0$. Note that this is provoked by the integer character of the decision variables because none of the conic constraints are binding at the optimal solution (see Figure~\ref{EjemploIlustrativo1Editado}, panel right).
According to expression (\ref{eq.pfdf}) the probabilities of infeasibility are $0.0404$ and $0.0019$, which are lower than the required probability $\Phi(-\beta_i)=0.1587$, i.e. the integer character of the variables involved provides a conservative solution for this particular case. We also check the probabilities of infeasibility by using a Monte Carlo sampling of 100000 realizations, obtaining estimated probabilities of infeasibility associated with reliability indexes equal to $1.7516$ and $2.8673$, respectively, which are very close to the exact values with relative errors of $0.32$\% and $0.68$\%, respectively.

In contrast, if the iterative method proposed in this paper is used, the convergence to the optimal solution within an pre-specified tolerance to $\varepsilon=10^{-6}$ is achieved in 2 iterations. The evolution of variables for the master and subproblems at every iteration is given in Table~\ref{TabIllEx}. Note that the algorithm converges to the same optimal values from problem (\ref{eq.fobjFIllusEx2})-(\ref{integralityIE2}).

\begin{table}
\begin{center}
  \begin{tabular}{cccccccc}
  \hline
  % after \\: \hline or \cline{col1-col2} \cline{col3-col4} ...
  $\nu$ & $\hat a_{11}$ & $\hat a_{12}$ & $\hat a_{21}$ & $\hat a_{22}$ & $\hat x_{1}$ & $\hat x_{2}$ \\
  \hline
1            &      1.00000 &      2.00000 &      2.00000 &      1.00000 &      1.00000 &      3.00000 \\
2            &      1.08496 &      2.19923 &      2.03780 &      1.07559 &      2.00000 &      2.00000 \\
  \hline
\end{tabular}
\end{center}
 \caption{Evolution of the iterative algorithm for the illustrative example.}\label{TabIllEx}
\end{table}

Note that in terms of computational time, results achieved using the QCP approach are slightly better. Nevertheless, we implemented all the problems using GAMS \citep{BrookeKMR:98} and did not make any special efforts to implement the individual steps of our algorithm efficiently, for instance, by taking advantage that the problems differ slightly on the values of the parameters, specially at the latest iterations. Note also that GAMS takes some time to build the models and this is done at every iteration. Such savings could potentially improve the running time of the algorithm, but not change the number of iterations required. To compare computational performance in a more meaningful way we present the following case study.

\section{Case study: Optimal Truss Design}\label{s.sheexam}
\begin{figure}[t]
\begin{center}
 \includegraphics*[width=0.95\textwidth]{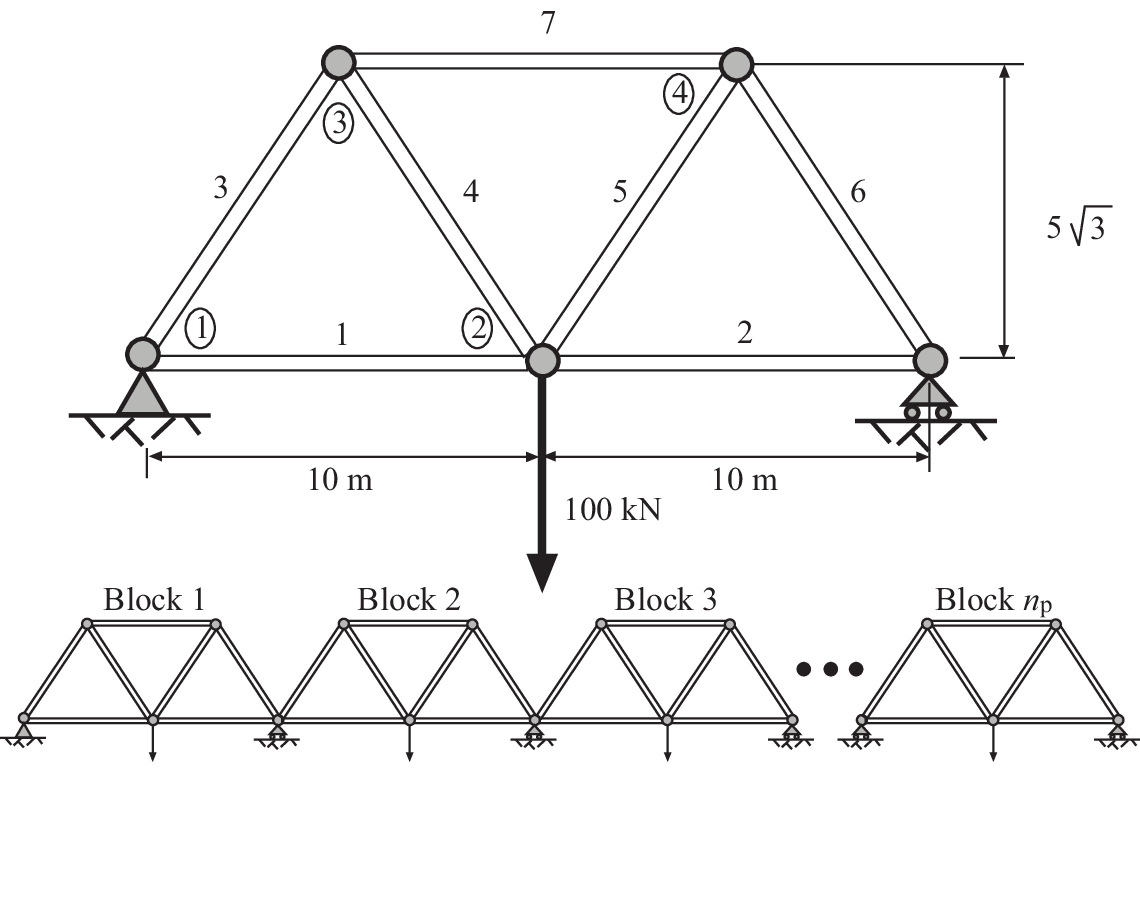}
\caption{\label{Truss} Truss design example.}
\end{center}
\end{figure}

This section considers an adapted example about optimal truss design previously used in different works \citep{RockafellarR:10,MinguezJTL:13,MinguezJTC:13}, a simple supported truss with 7 elements (bars) as shown in the upper part of Figure~\ref{Truss}. Yield stress of all members $\tilde a_i;i=1,\ldots,7$ are random variables with the following mean and variances: $\mbox{E}[\tilde{\bfg a}] = (-100,$ $-100,-200,$ $-200,-200,-200,-200)^T;N/mm^2$ and $\mbox{E}[(\tilde{\bfg a}-\mbox{E}[\tilde{\bfg a}])^2] = (15^2,$ $15^2,40^2,$ $40^2,40^2,40^2,$ $40^2)^T;N^2/mm^4$. Note that we use negative values for yield stress because it is more convenient from the formulation perspective. There is a vertical load applied on the structure which is also normally distributed with mean $p=100kN$ and standard deviation $\sigma_p=40kN$.

The aim of the design problem is to determine the cross-sectional areas of the bars $x_i; i=1,\ldots,7$, so that the probability of failure of each bar due to the uncertainty on yield stress and load is at most $0.001$. Instead of working with failure probabilities, and since random parameters are normally distributed, we use relationship (\ref{eq.pfdf}) to define the protection level of each bar $\beta_i=3.09; i=1,\ldots,7$. The advantage of this example is that the problem can be easily augmented in size by simply replicating the same block structure as shown in the lower part of Figure~\ref{Truss}.
%%%, so that the optimal solution of each block is equal to the optimal solution of one block, i.e. we know in advance the optimal solution of the problem if more than one block is used.
Assuming that there are $n_b$ blocks, the robust formulation of the design problem can be written as follows:
\begin{eqnarray}
&& \mm{Minimize}{x_{ik};\forall i;\forall k}{\displaystyle\sum_{k=1}^{n_b}\sum_{i=1}^7 c_ix_{ik}},\label{eq.fobjTruss}\\
% \nonumber to remove numbering (before each equation)
&&\hspace*{0.2cm}\mbox{subject to }p_k/\tau_i+a_{ik}x_{ik}  \le 0;(a_{ik},p_k)\in U_{ik};\forall i;\forall k\label{eq.inconsF1Truss}\\
&& \hspace*{0.2cm} x_{ik}\in\{0.5,0.6,0.7,\ldots,1.8,1.9,2\}; \forall i;\forall k,\label{eq.inconsF2Truss}
\end{eqnarray}
where $c_i=1$ are the cost coefficients, and $\tau_i$ are factors that depend on geometry and the load which are equal to $\tau_i=1/(2\sqrt{3})$ for $i=1,2$, and $\tau_i=1/(\sqrt{3})$ for $i=3,4,5,6,7$. The left hand side of constraints (\ref{eq.inconsF1Truss}) correspond to the difference between the actual stresses induced by the vertical load and the actual strength of the bars, note that the negative sign is implicitly included in the yield stress parameter. The optimal solution of one block in terms of decision variables is the same for all blocks, regardless of the number of blocks $n_b$ selected, for this reason we can use this example to compare computational performance between the traditional QCP and the proposed decomposition method on problems of different size, and checking afterwards if the optimal solution is attained. Note that that cross sectional areas can only take specific values from a given catalogue, i.e. $x_{ik}\in\{0.5,0.6,0.7,\ldots,1.8,1.9,2\}$.

The traditional formulation (\ref{eq.fobjF})-(\ref{eq.integralityF}) proposed by \cite{Ben-TalN:99} for this example becomes:
\begin{eqnarray}
&& \mm{Minimize}{x_{ik};\forall i;\forall k}{\displaystyle\sum_{k=1}^{n_b}\sum_{i=1}^7 c_ix_{ik}},\label{eq.fobjTrussRC}\\
% \nonumber to remove numbering (before each equation)
&&\hspace*{0.2cm}\mbox{subject to }p/\tau_i+a_{i}x_{ik}+\beta_i\sqrt{(\sigma_p/\tau_i)^2+\sigma_i^2 x_{ik}^2}\le0;\;\forall i;\forall k\label{eq.inconsF1TrussRC}\\
&&  \hspace*{0.2cm} x_{ik}\in\{0.5,0.6,0.7,\ldots,1.8,1.9,2\};\;\forall i;\forall k,\label{eq.inconsF2TrussRC}
\end{eqnarray}
where $p$ is the nominal value of all loads, $a_{i};i=1,\ldots,7$ and $\sigma_i;i=1,\ldots,7$ are, respectively, the nominal and standard deviation of yield stresses
associated with bars. Problem (\ref{eq.fobjTrussRC})-(\ref{eq.inconsF2TrussRC}) corresponds to a mixed-integer conic quadratic problem.

In contrast, the master and subproblems proposed in this paper for problem (\ref{eq.fobjTrussRC})-(\ref{eq.inconsF2TrussRC}) are defined as follows:
\begin{eqnarray}
&&\mm{Minimize}{x_{ik};\forall i;\forall k}{\displaystyle\sum_{k=1}^{n_b}\sum_{i=1}^7 c_ix_{ik}},\label{eq.fobjTrussMaster}\\
% \nonumber to remove numbering (before each equation)
&&\hspace*{0.2cm}\mbox{subject to }\hat p_{ik}/\tau_i+\hat a_{ik}x_{ik} \le 0;\;\forall i;\forall k\label{eq.inconsF1TrussMaster}\\
&& \hspace*{0.2cm}x_{ik}\in\{0.5,0.6,0.7,\ldots,1.8,1.9,2\}; \forall i;\forall k,\label{eq.inconsF2TrussMaster}
\end{eqnarray}
and
\begin{equation}\label{suproTruss}
\left.
\begin{array}{l}
\mm{Maximize}{\tilde p_k,\tilde a_{ik}}{} \tilde p_k/\tau_i+\tilde a_{ik}x_{ik}\\
\hspace*{0.2cm}\mbox{subject to}   \\
% \nonumber to remove numbering (before each equation)
 %%% z_{pk} = & \fraca{\tilde p_k-p}{\sigma_p}\\
%%% z_{ik} = & \fraca{\tilde a_{ik}-a_i}{\sigma_i}\\
 \hspace*{0.2cm}    \left(\fraca{\tilde p_k-p}{\sigma_p}\right)^2+\left(\fraca{\tilde a_{ik}-a_i}{\sigma_i}\right)^2 \le  \beta_i^2
\end{array}\right\}\forall i\forall k.
\end{equation}

The optimal solution of subproblem (\ref{suproTruss}) according to (\ref{eq.ahati}) is:
\begin{eqnarray}
    \left(\begin{array}{c}\hat p_k\\ \hat a_{ik}\\\end{array}\right)&=&\left(\begin{array}{c} p\\ a_{i}\\\end{array}\right)+\beta_i \fraca{(\sigma_p^2/\tau_i+x_{ik} \sigma_i^2)}{\sqrt{(\sigma_p^2/\tau_i^2+x_{ik}^2 \sigma_i^2)}};\;\forall i\forall k.\label{eq.finTrussMaster}
\end{eqnarray}

The optimal solution associated with the RC and MIQCP problem (\ref{eq.fobjTrussRC})-(\ref{eq.inconsF2TrussRC}) depends on the number of blocks considered $n_b$ and it is equal to:
\begin{equation}\label{QCPsol}
\begin{array}{rl}
  \hat f = & 7.3 n_b \\
  \hat{\bfg x}^{(k)}= & (0.9,\;0.9,\;1.1,\;1.1,\;1.1,\;1.1,\;1.1)^T;\;\forall k=1,\ldots,n_b.
\end{array}
\end{equation}

The corresponding problem has $7 \times n_b$ quadratic constraints, and $7 \times n_b$ integer variables with 15 different possible values and one continuous variable.

In order to compare the above procedures, the following problems have been solved using different current mathematical modelling solvers and different number of blocks $n_b\in\{10,100,$ $1000,10000,100000\}$:
\begin{enumerate}
  %%%  \item QCP problem (\ref{eq.fobjTrussRC})-(\ref{eq.inconsF2TrussRC}) using CONOPT \citep{Drud:96} and COUENNE (http://www.coin-or.org/) solvers.
%%%    \item LP$^{(1)}$ problem (\ref{eq.fobjTrussMaster})-(\ref{eq.finTrussMaster}) using the decomposition procedure proposed in this paper. Note that master LP problems are solved using CBC (https://projects.coin-or.org/Cbc), CONOPT, CPLEX (http://www.ilog.com) and MINOS \citep{MurtaghS:98} solvers.
%%%
%%%     \item LP$^{(2)}$ problem (\ref{eq.fobjBS04})-(\ref{eq.inconsF6BS04}) using the robust formulation proposed by  \citep{BertsimasS:04}. Note that master LP problems are solved using CBC (https://projects.coin-or.org/Cbc), CONOPT, CPLEX (http://www.ilog.com) and MINOS \citep{MurtaghS:98} solvers.

    \item MIQCP problem (\ref{eq.fobjTrussRC})-(\ref{eq.inconsF2TrussRC}) using BONMIN (COIN-OR Bonmin 24.2.3 r46072, \citep{Bonamietal:08}), COUENNE (COIN-OR Couenne  24.2.3 r46072), DICOPT (CONOPT3 version 3.15P and CPLEX 12.6, \citep{Rosenthal:08}), SBB (CONOPT3 version 3.15P, \citep{Rosenthal:08}) and CPLEX 12.6 solvers.

    \item MIP problem (\ref{eq.fobjTrussMaster})-(\ref{eq.finTrussMaster}) using the decomposition procedure proposed in this paper. Note that the master MIP problems are solved using BONMIN (COIN-OR Bonmin 24.2.3 r46072) and CPLEX 12.6 solvers.

        %%% \item MIP$^{(2)}$ problem (\ref{eq.fobjBS04})-(\ref{eq.inconsF6BS04})  and using the robust formulation proposed by  \citep{BertsimasS:04}. Note that the master MIP problems are solved using BONMIN and CPLEX solvers.
\end{enumerate}

All computations have been performed on an Intel Xeon E7-4820 computer with
four processors clocking at 2GHz and 756GB of RAM under GAMS release 24.2.3. It is worth mentioning that all results associated with the decomposition procedure (MIP) are obtained
after three iterations of the proposed method, using a feasibility tolerance of $\epsilon=10^{-6}$. We imposed a time limit of 9600 seconds so that if the solver does not find an optimal solution within that time window the process is stopped.

\begin{table}
\begin{center}
\renewcommand{\arraystretch}{1.0}
  \begin{tabular}{ccccccc}
  \cline{3-7}
&  &\multicolumn{5}{c}{CPU Time (seconds)}\\
  \cline{2-7}
Problem   & Solver & $n_b=10$ & $n_b=10^2$ &$n_b=10^3$ & $n_b=10^4$ & $n_b=10^5$\\
%%%      \hline
%%%      \hline
%%%(\ref{eq.fobjTrussRC})-(\ref{eq.inconsF2TrussRC})     & CONOPT & 0.719 & 0.703 & 51.203 & 9054.015 \\
%%%QCP                                                  & COUENNE & 14.032 & 1832.890 & 9542.219 & 9600$^{\ast}$ \\
%%%     \hline
%%%   & CBC & 0.469 & 0.500 & 1.828 & 68.312 \\
%%%(\ref{eq.fobjTrussMaster})-(\ref{eq.finTrussMaster})   & CONOPT & 0.578  & 0.609 & 1.719 & 67.656 \\
%%% LP$^{(1)}$  & CPLEX & 0.500 & 0.484 & 1.703 &  67.313 \\
%%%   & MINOS & 0.531 & 0.500 & 2.703 & 195.141 \\
%%%   \hline
%%%      & CBC & 0.110 & 0.187 & 2.718 & 134.313\\
%%%(\ref{eq.fobjBS04})-(\ref{eq.inconsF6BS04})  & CONOPT & 0.125 & 2.922 & 286.828 & 9600$^{\ast}$ \\
%%% LP$^{(2)}$  & CPLEX & 0.125 & 0.265 & 1.297 & 31.313\\
%%%   & MINOS & 0.110 & 0.453 & 46.468 & 7888.547 \\
   \hline
   \hline
     & BONMIN & 0.78$^\S$ & 1.11$^\S$ & 4.57$^\S$ & 45.81$^\S$ & 675.46$^\S$\\
 (\ref{eq.fobjTrussRC})-(\ref{eq.inconsF2TrussRC}) & COUENNE & 0.42 & 1.34 & 67.40 & 9600$^{\ast}$ & 9600$^{\ast}$\\
 MIQCP  & DICOPT & 1.49  & 50.581 & 9517.75 & 9585.99 & 9600$^{\ast}$\\
   & SBB & 0.85 & 10.76 & 87.52 & 1789.32  & 9600$^{\ast}$\\
   & CPLEX  & 0.43 & 0.48 & 1.34 & 22.52 &  1113.37\\
   \hline
(\ref{eq.fobjTrussMaster})-(\ref{eq.finTrussMaster})     & BONMIN & 1.13 & 1.18 & 2.59 & 125.437 & 4356.91\\
 MIP  & CPLEX & 1.33 & 1.34 & 2.23  & 105.71 & 4486.17\\
%%%  \hline
%%% (\ref{eq.fobjBS04})-(\ref{eq.inconsF6BS04})      & BONMIN & 0.703 &  2.156 & 42.969 & 84.781$^\S$\\
%%% MIP$^{(2)}$  & CPLEX & 0.234 & 0.344 & 4.172 & 363.500 \\
  \hline \hline
  %%%\multicolumn{7}{l}{{\footnotesize $(1) \Longrightarrow$ Proposed iterative method.}}\\
  %%%\multicolumn{6}{l}{{\footnotesize $(2) \Longrightarrow$ Bertsimas and Sim (2004) method \citep{BertsimasS:04}.}}\\
  \multicolumn{7}{l}{{\footnotesize $\ast:$ Maximum time limit reached and no optimal solution found.}}\\
  \multicolumn{7}{l}{{\footnotesize $\S:$ No optimal solution found within maximum time limit.}}\\
  \hline
  \end{tabular}
  \caption{Computational results of the case study using different methods, solvers and problem types and sizes.}\label{t.res1}
\end{center}
\end{table}

Table~\ref{t.res1} provides the computational times in seconds taken for each solver to reach the optimal solution for the different cases and problems considered.
According to these
results, the following observations are pertinent:
\begin{enumerate}
 %%%\item Computational time for QCP formulation increases faster with respect to the size of the problem than LP$^{(1)}$ formulation proposed in this paper. Note in the first two rows of Table~\ref{t.res1} that both solvers almost reached the time limit of 9600 seconds, in fact, solver COUENNE reached the time limit of 9600 without achieving the optimal solution of the problem.
%%%\item Computational time for LP$^{(1)}$ formulation associated with the proposed iterative method is considerably faster than the equivalent QCP formulation, as shown in rows 3 to 6 of Table~\ref{t.res1}. Note that all solvers, even those for non-linear programming, reached the solution in a reasonable amount of time.
%%%    \item Computational time for LP$^{(2)}$ formulation associated with the method proposed by \citep{BertsimasS:04} is analogous to that obtained for LP$^{(1)}$ if linear programming solvers CBC and CPLEX are used. For non-linear programming solvers the computational time grows rapidly, in fact solver CONOPT reached the time limit of 9600 without achieving the optimal solution of the problem.
%%%  \item When the number of block structures is relatively low, i.e. below 100, computational times related to the proposed decomposition procedure might be higher than the alternative quadratic formulations. This effect is produced by the deficiencies in the implementation of the algorithm, which includes model generation times and other times that could be avoided if an ad-hoc solver implementing the proposed procedure is programmed.

\item Computational time for MIQCP formulation, analogously to the QCP formulation,  increases exponentially with respect to the size of the problem. For 100000 blocks only CPLEX solver succeeds on finding the optimal solution within the maximum time frame of 9600 seconds considered in this work, this result confirms conclusion by \citep{BertsimasS:04} that robust optimization using ellipsoidal uncertainty sets is not particularly attractive if integer variables are involved for most of the solvers.

\item Computational time for MIP formulation associated with the proposed iterative method allows solving robust optimization problems using ellipsoidal uncertainty sets provided that the appropriate mixed-integer solver, such as BONMIN or CPLEX, is used. Both solvers perform similarly for this particular example. Note that computational times are considerably lower than those related to MIQCP formulation except for CPLEX solver.

 %%%\item Computational time for MIP$^{(2)}$ formulation associated with the method proposed by \citep{BertsimasS:04} is analogous to
%%%      that obtained for MIP$^{(1)}$ provided that the appropriate mixed-integer solver, such as CPLEX, is used. Note that when the size of the problem grows, MIP$^{(2)}$ formulation seems to be a bit slower than MIP$^{(1)}$, although this result is not statistically significant. However, we can conclude that both approaches MIP$^{(1)}$ and MIP$^{(2)}$ are alternative formulations to solve robust optimization problems involving mixed-integer variables.

\end{enumerate}

%%%To further reinforce our conclusions, specially the last one in terms of tractability, we run an example using $n_b=30000$ blocks and solved the robust problems associated with MIP$^{(1)}$ and MIP$^{(2)}$ formulations. Problem MIP$^{(1)}$ at the last iteration has $630,001$ equations, $210001$ continuous variables and $210,000$ discrete variables, and it was solved in $1869.609$ seconds of CPU time. In contrast, problem MIP$^{(2)}$ has $1050001$ equations, $1050001$ continuous variables and $210000$ discrete variables, and it was solved in $1886.015$ seconds of CPU time. Note that both computational times are of the same order of magnitude.
%%%
%%%Finally, we run an additional example using $n_b=10^5$ blocks and solved the robust problems associated with MIP$^{(1)}$ and MIP$^{(2)}$ formulations. Problem MIP$^{(1)}$ at the last iteration has $2,100,001$ equations, $700,001$ continuous variables and $700000$ discrete variables, and it was solved in $27675.953$ seconds of CPU time. However, problem MIP$^{(2)}$, which has $3500001$ equations, $3500001$ continuous variables and $700000$ discrete variables, could not be solved due to memory problems. CPLEX solver produced an error 1001: Out of memory.

Note that the proposed algorithm is competitive with respect most of the solvers available in the state-of-the-art except for CPLEX 12.6 solver, that can now handle mixed-integer second-order cone programs. In this particular case, our method is worse in terms of computing time. However, as mentioned in the illustrative example, we implemented all the problems using GAMS \citep{BrookeKMR:98} and did not make any special efforts to implement the individual steps of our algorithm efficiently. Nevertheless, in work \citep{BertsimasDL:16} an intensive computational study is made comparing a modified version of the method proposed in this paper to improve computational efficiency, and concluded that there is no dominant method when dealing with robust mixed-integer problems, which make this cutting plane methods a plausible alternative for solving this kind of problems.

\section{Discussion of possible extensions}\label{discussion}
Although the problem dealt with in this paper is useful for many different applications, this type of robust optimization models is known to be conservative. Therefore, an interesting feature for further research is the use of joint probabilistic constraints, where restrictions (\ref{eq.inconsRCi2}) are replaced by constraint:
\begin{equation}\label{eq.inconsRCi2b}
    \mbox{Prob}\left( {\tilde{\bfg a}}^{(i)^T} {\bfg x}  \le  0;\;\forall i=1,\ldots,m\right)\ge \Phi(\beta).
\end{equation}

Combining the methods presented in \citep{AckooijHMZ:11} and \citep{Bremer:15}, the model including this alternative constraint can likely be fully solved under the assumption that the coefficients of ${\bfg A}$ follow a multivariate distribution function. The difference between both approaches is apparent if we consider the truss case study. In this example, we optimize cross sectional areas assuming that the probability of failure of each bar must be lower or at least equal to the target probability of 0.001. However, the collapse of the block structure might occur if any of the bars fail, i.e. we are dealing with a series structural system, so the probability of collapse is greater than 0.001. For instance, assuming that the bar strengths are independent, the probability of collapse is equal to $P_{collapse}=1-(1-0.001)^7= 0.00698$. Therefore, it is more convenient to optimize the structure using the joint chance constraint (\ref{eq.inconsRCi2b}), which represents the probability of survival of the block structure as a whole. An additional advantage of using this alternative is that it is possible to consider correlations among bar strengths, which is more realistic.

\section{Conclusions}\label{Conclu}
Based on decomposition
techniques, this paper proposes an iterative method for solving RC of uncertain mixed-integer linear programs with ellipsoidal uncertainty sets. The method
is specially suitable for problems where first and second order moments of the probability distributions of the uncertain parameters involved are available. In addition, the proof of convergence and expressions for the probability of constraint violation are given, which allows expressing the robust counterpart problem as a chance constraint mathematical programming problem.

Although last versions of state-of-the-art solvers, such as CPLEX 12.6, can now handle mixed-integer second-order cone programs efficiently, as shown in the computational study, the method proposed in this paper is also robust and efficient and can be considered an alternative for solving this kind of problems. Besides, it is demonstrated in the current literature that slight modifications and an ad-hoc implementation of the algorithm proposed in this work make both ways of solving these problems analogous in terms of computing time.

%%%In summary, the method proposed in this paper allows the applicability of the RC developed by \citep{Ben-TalN:99} for a wider range of problems, even making it attractive for solving robust discrete optimization models. This decomposition procedure makes RC with ellipsoidal uncertainty to be an alternative method in terms of tractability with respect to the one proposed by \cite{BertsimasS:04}, which uses cardinality constrained uncertainty sets.
%%%
%%%Further research includes the extension of the proposed decomposition technique to linear programming formulations, and its comparison in terms of performance with respect to interior-point algorithms.

\section*{Acknowledgments}
Dr. Casero-Alonso has been sponsored by Ministerio de Econom\'{\i}a y Competitividad and grant contract FEDER MTM2013-47879-C2-1-P. The authors want to express their gratitude to Dr. Pablo Pedregal for his comments and discussions about the paper, which certainly improve the original manuscript.
The authors also thank the referees for their comments, suggestions and corrections, which contributed to enhance the quality of the paper.

\bibliographystyle{elsarticle-num}
%%%\bibliography{BiblioIngAll}

\begin{thebibliography}{10}
\expandafter\ifx\csname url\endcsname\relax
  \def\url#1{\texttt{#1}}\fi
\expandafter\ifx\csname urlprefix\endcsname\relax\def\urlprefix{URL }\fi
\expandafter\ifx\csname href\endcsname\relax
  \def\href#1#2{#2} \def\path#1{#1}\fi

\bibitem{Soyster:73}
A.~L. Soyster, Convex programming with set-inclusive constraints and
  applications to inexact linear programming, Oper. Res. 21~(5) (1973)
  1154--1157.

\bibitem{El-GhaouiL:97}
L.~El-Ghaoui, H.~Lebret, Robust solutions to least-squares problems with
  uncertain data, SIAM J. Matrix Anal. Appl. 18~(4) (1997) 1035--1064.

\bibitem{El-GhaouiOL:98}
L.~El-Ghaoui, F.~Oustry, H.~Lebret, Robust solutions to uncertain semidefinite
  programs, SIAM J. on Optimization 9~(1) (1998) 33--52.

\bibitem{Ben-TalN:98}
A.~Ben-Tal, A.~Nemirovski, Robust convex optimization, Mathematics of
  Operations Research 23~(4) (1998) 769--805.

\bibitem{Ben-TalN:99}
A.~Ben-Tal, A.~Nemirovski, Robust solutions of uncertain linear programs,
  Operations Research Letters 25~(1) (1999) 1 -- 13.

\bibitem{Ben-TalN:00}
A.~Ben-Tal, A.~Nemirovski, Robust solutions of linear programming problems
  contaminated with uncertain data, Mathematical Programming 88~(3) (2000)
  411--424.

\bibitem{BertsimasS:04}
D.~Bertsimas, M.~Sim, The price of robustness, Oper. Res. 52~(1) (2004) 35--53.

\bibitem{BirgeL:97}
J.~R. Birge, F.~Louveaux, Introduction to Stochastic Programming, Springer
  Verlag, New York, 1997.

\bibitem{RockafellarU:00}
R.~T. Rockafellar, S.~Uryasev, Optimization of conditional value-at risk, J.
  Risk 2~(3) (2000) 21--41.

\bibitem{RockafellarU:02}
R.~T. Rockafellar, S.~Uryasev, Conditional value-at risk for general loss
  distributions, J. Bank. Finan. 26~(7) (2002) 1443--1471.

\bibitem{GabrelMT:14}
V.~Gabrel, C.~Murat, A.~Thiele, Recent advances in robust optimization: {A}n
  overview, European Journal of Operational Research 235 (2014) 471--483.

\bibitem{HoushOS:11}
M.~Housh, A.~Ostfeld, U.~Shamir, Optimal multiyear management of a water supply
  system under uncertainty: Robust counterpart approach, Water Resources
  Research 47~(10) (2011)  n/a--n/a.
\newblock http://dx.doi.org/10.1029/2011WR010596

\bibitem{PerelmanHO:13}
L.~Perelman, M.~Housh, A.~Ostfeld, Robust optimization for water distribution
  systems least cost design, Water Resources Research 49~(10) (2013)
  6795--6809.

\bibitem{Frangopol:95}
D.~M. Frangopol, Reliability-based optimum structural design, in: C.~Sundarajan
  (Ed.), Probabilistic Structural Mechanics Handbook, Chapmam \&\ Hall, New
  York, 1995, Ch.~16, pp. 352--387.

\bibitem{Melchers:99}
R.~E. Melchers, Structural reliability analysis and prediction, 2nd Edition,
  John Wiley \&\ Sons, New York, 1999.

\bibitem{RoysetDP:01}
J.~O. Royset, A.~Der~Kiureghian, E.~Polak, Reliability-based optimal design of
  series structural systems, Journal of Engineering Mechanics, ASCE 127~(6)
  (2001) 607--614.

\bibitem{RoysetDP:01b}
J.~O. Royset, A.~Der~Kiureghian, E.~Polak, Reliability-based optimal structural
  design by the decoupling approach, Reliab. Eng. Syst. Saf. 73~(3) (2001)
  213–--221.

\bibitem{RoysetDP:06}
J.~O. Royset, A.~Der~Kiureghian, E.~Polak, Optimal design with probabilistic
  objective and constraints, Journal of Engineering Mechanics, ASCE 132~(1)
  (2006) 107--118.

\bibitem{HasoferL:74}
A.~M. Hasofer, N.~C. Lind, Exact and invariant second moment code format, J.
  Engrg. Mech. 100~(EM1) (1974) 111--121.

\bibitem{MinguezCG:11}
R.~M\'{\i}nguez, A.~J. Conejo, R.~Garc\'{\i}a-Bertrand, Reliability and
  decomposition techniques to solve certain class of stochastic programming
  problems, Reliability Engineering \& System Safety 96 (2011) 314--323.

\bibitem{RockafellarR:10}
R.~T. Rockafellar, J.~O. Royset, On buffered failure probability in design and
  optimization of structures, Reliability Engineering \& System Safety 95
  (2010) 499--510.

\bibitem{Floudas:95}
C.~A. Floudas, Nonlinear and Mixed-Integer Optimization. Fundamentals and
  Applications, Oxford University Press, New York, 1995.

\bibitem{ConejoCMG:06}
A.~J. Conejo, E.~Castillo, R.~M\'{\i}nguez, R.~Garc\'{\i}a-Bertrand,
  Decomposition techniques in mathematical programming. Engineering and science
  applications, Springer-Verlag Berlin Heidelberg, New York, 2006.

\bibitem{FischettiM:12}
M.~Fischetti, M.~Monaci, Cutting plane versus compact formulations for
  uncertain (integer) linear programs, Math. Program. Comput. 4~(3) (2012)
  239–-273.

\bibitem{BertsimasDL:16}
D.~Bertsimas, I.~Dunning, M.~Lubin, Reformulation versus cutting-planes for
  robust optimization: {A} computational study, Comput. Manag. Sci. 13~(2)
  (2016) 195–-217.

\bibitem{AckooijFD:16}
W.~van Ackooij, A.~Frangioni, W.~de~Oliveira, Inexact stabilized benders'
  decomposition approaches with application to chance-constrained problems with
  finite support, Computational Optimization and Applications 65~(3) (2016)
  637–-669.

\bibitem{Ben-TalEN:09}
A.~Ben-Tal, L.~El~Ghaoui, A.~Nemirovski, Robust Optimization, Princeton Series
  in Applied Mathematics, Princeton University Press, Princeton, 2009.

\bibitem{DelageY:10}
E.~Delage, Y.~Ye, Distributionally robust optimization under moment uncertainty
  with application to data-driven problems, Operations Research 58~(3) (2010)
  595--612.

\bibitem{BertsimasBC:11}
D.~Bertsimas, D.~B. Brown, C.~Caramanis, Theory and applications of robust
  optimization, SIAM Review 53~(3) (2011) 464–-501.

\bibitem{Pinar:00}
M.~Pinar, A simple duality proof in convex quadratic programming with a
  quadratic constraint, and some applications, European Journal of Operational
  Research 124~(1) (2000) 151--158.

\bibitem{Ackooij:15}
W.~van Ackooij, A comparison of four approaches from stochastic programming for
  large-scale unit-commitment, EURO Journal on Computational Optimization
  (2015) 1--29.

\bibitem{FletcherL:94}
R.~Fletcher, S.~Leyffer, Solving mixed integer nonlinear programs by outer
  approximation, Mathematical Programming 66~(3) (1994) 327--349.

\bibitem{DuranG:86}
M.~Duran, I.~E. Grossmann, An outer-approximation algorithm for a class of
  mixed integer nonlinear programs, Mathematical Programming 36 (1986)
  307--339.

\bibitem{WesterlundP:95}
T.~Westerlund, F.~Pettersson, An extended cutting plane method for solving
  convex {MINLP} problems, Computers \& Chemical Engineering 19 (1995)
  131--136.

\bibitem{Freudenthal:56}
A.~M. Freudenthal, Safety and the probability of structural failure,
  Transactions, ASCE 121 (1956) 1337--1397.

\bibitem{RackwitzF:78}
R.~Rackwitz, B.~Fiessler, Structural reliability under combined load sequences,
  Comput. Struct. 9 (1978) 489--494.

\bibitem{Ditlevsen:81}
O.~Ditlevsen, Principle of normal tail approximation, J. Engineering Mechanics
  Div., ASCE 107~(6) (1981) 1191--1208.

\bibitem{HohenbichlerR:81}
M.~Hohenbichler, R.~Rackwitz, Non-normal dependent vectors in structural
  safety, J. Engineering Mechanics Div., ASCE 107~(6) (1981) 1227--1238.

\bibitem{MinguezC:09}
R.~M\'{\i}nguez, M.~Castillo, Reliability-based optimization in engineering
  using decomposition techniques and {FORMS}, Structural Safety 31~(3) (2009)
  214--223.

\bibitem{BrookeKMR:98}
A.~Brooke, D.~Kendrick, A.~Meeraus, R.~Raman, GAMS: A user's guide, GAMS
  Development Corporation, Washington, 1998.

\bibitem{MinguezJTL:13}
R.~M\'{\i}nguez, F.~F. Jaime, A.~Tom\'as, J.~L. Lara, Iterative scenario
  reduction technique to solve reliability-based optimization problems using
  the buffered failure probability, in: Proceedings of the 11th International
  Conference on Structural Safety \& Reliability ({ICOSSAR 2013}), New York,
  NY, 2013.

\bibitem{MinguezJTC:13}
R.~M\'{\i}nguez, F.~F. Jaime, A.~Tom\'as, E.~Castillo, New insights on the
  buffered failure probability risk measure for optimal structural design, in:
  Proceedings of the 11th International Conference on Structural Safety \&
  Reliability ({ICOSSAR 2013}), New York, NY, 2013.

\bibitem{Bonamietal:08}
P.~Bonami, L.~Biegler, A.~Conn, G.~Cornuejols, I.~Grossmann, C.~Laird, J.~Lee,
  A.~Lodi, F.~Margot, N.~Sawaya, A.~Waechter, An algorithmic framework for
  convex mixed integer nonlinear programs, Discrete Optimization 5~(2) (2008)
  186--204.

\bibitem{Rosenthal:08}
R.~E. Rosenthal, GAMS: A user's guide, GAMS Development Corporation,
  Washington, 2008.

\bibitem{AckooijHMZ:11}
W.~van Ackooij, R.~Henrion, A.~M\"oller, R.~Zorgati, On joint probabilistic
  constraints with gaussian coefficient matrix, Operations Research Letters 39
  (2011) 99--102.

\bibitem{Bremer:15}
I.~Bremer, R.~Henrion, A.~M\"oller, Probabilistic constraints via {SQP} solver:
  {A}pplication to renewable energy management problem, Computational
  Management Science 12 (2015) 435--459.

\end{thebibliography}

%%%\newpage
%%%
%%%\begin{figure}[t]
%%%\begin{center}
%%%%\resizebox{6.0cm}{!}{\includegraphics*{Transdens.eps}}
%%%\includegraphics*[width=0.8\textwidth]{Transdens.eps}
%%%\caption{\label{failureregion} Joint probability distribution
%%%function contours, limit-state equation, and design points in: a)
%%%the initial random space $\tilde{\bfg a}^{(i)}$, and b) the unit standard
%%%normal random space ${\bfg z}^{(i)}$.}
%%%\end{center}
%%%\end{figure}
%%%%
%%%%\newpage
%%%%
%%%%\begin{figure}[t]
%%%%\begin{center}
%%%% \includegraphics*[width=0.95\textwidth]{EjemploIlustrativoBinario}
%%%%\caption{\label{EjemploIlustrativo1Editado} Graphical illustration of illustrative example: panel on the left shows random instances of uncertain constraints, while panels on the right show the conic robust contraints that define the feeasibility region and the optimal solution.}
%%%%\end{center}
%%%%\end{figure}
%%%%
%%%%\newpage
%%%%
%%%%\begin{figure}[t]
%%%%\begin{center}
%%%% \includegraphics*[width=0.95\textwidth]{Truss.eps}
%%%%\caption{\label{Truss} Truss design example.}
%%%%\end{center}
%%%%\end{figure}
%%%%

\end{document}